\documentclass[12pt,oneside,reqno]{amsart}
\usepackage{graphicx}
\usepackage{mathrsfs}
\usepackage{stmaryrd}
\usepackage{amsfonts}
\usepackage{enumerate,amsmath,amssymb,amsthm}
\newcommand{\dif}{\mathrm{d}}

\newcommand{\be}{\begin{eqnarray}}
\newcommand{\ee}{\end{eqnarray}}
\newcommand{\ce}{\begin{eqnarray*}}
\newcommand{\de}{\end{eqnarray*}}
\newtheorem{theorem}{Theorem}[section]
\newtheorem{lemma}[theorem]{Lemma}
\newtheorem{remark}[theorem]{Remark}
\newtheorem{definition}[theorem]{Definition}
\newtheorem{proposition}[theorem]{Proposition}
\newtheorem{Example}[theorem]{Example}
\newtheorem{corollary}[theorem]{Corollary}

\def\t{\theta}
\def\a{\alpha}

\def\[{{\Big[}}
\def\]{{\Big]}}
\def\<{{\langle}}
\def\>{{\rangle}}
\def\({{\Big(}}
\def\){{\Big)}}

\def\tr{{\rm tr}}

\def\no{\nonumber}
\def\bt{\begin{theorem}}
\def\et{\end{theorem}}
\def\bl{\begin{lemma}}
\def\el{\end{lemma}}
\def\br{\begin{remark}}
\def\er{\end{remark}}
\def\bx{\begin{Example}}
\def\ex{\end{Example}}
\def\bd{\begin{definition}}
\def\ed{\end{definition}}
\def\bp{\begin{proposition}}
\def\ep{\end{proposition}}
\def\bc{\begin{corollary}}
\def\ec{\end{corollary}}
\def\cA{{\mathcal A}}

\def\cG{{\mathcal G}}
\def\cH{{\mathcal H}}

\def\cM{{\mathcal M}}
\def\cN{{\mathcal N}}

\def\cQ{{\mathcal Q}}

\def\cU{{\mathcal U}}
\def\cV{{\mathcal V}}

\def\mE{{\mathbb E}}

\def\mN{{\mathbb N}}

\def\mP{{\mathbb P}}

\def\mR{{\mathbb R}}
\def\mS{{\mathbb S}}

\def\sB{{\mathscr B}}

\def\sF{{\mathscr F}}

\def\geq{\geqslant}
\def\leq{\leqslant}

\begin{document}

\allowdisplaybreaks

\title{Path independence of the additive functionals for stochastic differential equations driven by G-L\'evy processes*}

\author{Huijie Qiao$^{1,2}$ and Jiang-Lun Wu$^3$}



\thanks{*This work was partly supported by NSF of China (No. 11001051, 11371352, 11671083) and China Scholarship Council under Grant No. 201906095034.}


\date{}

\dedicatory{1. School of Mathematics,
Southeast University\\
Nanjing, Jiangsu 211189,  China\\
2. Department of Mathematics, University of Illinois at
Urbana-Champaign\\
Urbana, IL 61801, USA\\
hjqiaogean@seu.edu.cn\\
3. Department of Mathematics, Computational Foundry, Swansea University\\
Bay Campus, Swansea SA1 8EN, UK\\
j.l.wu@swansea.ac.uk}

\begin{abstract} 
In the paper, we consider a type of stochastic differential equations driven by G-L\'evy processes. We prove that a kind of their additive functionals has path independence and extend some known results.
\end{abstract}

\maketitle \rm

\section{Introduction}

Recently, development of mathematical finance forces the appearance of a type of processes-G-Brownian motions(\cite{p2}). And then the related theory, such as stochastic 
calculus and stochastic differential equations (SDEs in short) driven by G-Brownian motions, are widely studied(\cite{g, p2, peng, q1}). However, in some financial models, volatility uncertainty makes G-Brownian motions insufficient 
for simulating these models. One important reason lies in the continuity of their paths with respect to the time variable. So, Hu-Peng \cite{hp} solved the problem by introducing G-L\'evy processes. 
And the type of processes has discontinuous (right continuous with left limits) paths. Later, Paczka \cite{pk1} defined the It\^o-L\'evy stochastic integrals, deduced the It\^o formula, established 
SDEs driven by G-L\'evy processes and stated the existence and uniqueness of solutions for these equations under lipschitz conditions. Most recently, under non-lipschitz conditions, Wang-Gao \cite{wg} 
proved well-definedness of SDEs driven by G-L\'evy processes and investigated exponential stability of their solutions. Here, we follows up the line in \cite{wg}, define the additive functionals of SDEs driven by G-L\'evy processes and study their path independence.

Concretely speaking, we consider the following SDEs on $\mR^d$:
\be
\dif Y_t=b(t,Y_t)\dif t +h_{ij}(t,Y_t)\dif \<B^i, B^j\>_t+\sigma(t,Y_t)\dif B_t+\int_{\mR^d\setminus\{0\}}f(t,Y_t,u)L(\dif t,\dif u),
\label{glesde}
\ee
where $B$ is a G-Brownian motion, $\<B^i, B^j\>_t$ is the mutual variation process of $B^i$ and $B^j$ for $i,j=1,2,\cdots,d$ and $L(\dif t,\dif u)$ is a G-random measure (See Subsection \ref{itointe}). 
The coefficients $b: [0,T]\times\mR^d\mapsto\mR^d$, $h_{ij}=h_{ji}: [0,T]\times\mR^d\mapsto\mR^d$, $\sigma: [0,T]\times\mR^d\mapsto\mR^{d\times d}$ and $f: [0,T]\times\mR^d\times(\mR^d\setminus\{0\})\mapsto\mR^d$ 
are Borel measurable. Here and hereafter we use the convention that the repeated indices stand for the summation. Thus, under (${\bf H}^1_{b,h,\sigma,f}$)-(${\bf H}^2_{b,h,\sigma,f}$) in Subsection \ref{sdegle}, 
by \cite[Theorem 3.1]{wg}, we know that Eq.(\ref{glesde}) has a unique solution $Y_t$. And then we introduce the additive functionals of $Y_t$ and define path independence of these functionals. Finally, we prove 
that these functionals have path independence under some assumption.

Next, we say our motivations. First, we mention that Ren-Yang \cite{ry} proved path independence of additive functionals for SDEs driven by G-Brownian motions. Since these equations can not satisfy 
the actual demand very well, to extend them becomes one of our motivations. Second, by analyzing some special cases, we surprisingly find that, we can express explicitly these additive functionals. 
However, in our known results (\cite{qw1,qw2,qw3}), it is difficult to express explicitly additive functionals. Therefore, this is the other of our motivations.

This paper is arranged as follows. In Section \ref{pre}, we introduce G-L\'evy processes, the It\^o-L\'evy stochastic integrals, SDEs driven by G-L\'evy processes, additive functionals, path independence and some related results. The main results and their proofs are placed in Section \ref{main}. Moreover, we analysis some special cases and compare our result with some known results (\cite{qw1,qw2,qw3, ry}) in Subsection \ref{com}. 
 
\section{Preliminary}\label{pre}

In the section, we introduce some concepts and results used in the sequel.

\subsection{Notation}

In the subsection, we introduce notations used in the sequel. 

For convenience, we shall use $\mid\cdot\mid$ and $\parallel\cdot\parallel$  for norms of vectors and matrices, respectively. Furthermore, let $\langle\cdot$ , $\cdot\rangle$ denote the scalar product in $\mR^d$. Let $Q^*$ denote the transpose of the matrix $Q$.

Let $lip(\mR^n)$ be the set of all Lipschitz continuous functions on $\mR^n$ and $C_{b, lip}(\mR^d)$ be the collection of all bounded and Lipschitz continuous functions on $\mR^d$. Let $C_b^3(\mR^d)$ be the space of bounded and three times continuously differentiable functions with bounded derivatives of all orders less than or equal to $3$.

\subsection{G-L\'evy processes}\label{glevy}

In the subsection, we introduce G-L\'evy processes.(c.f.\cite{hp})

Let $\Omega$ be a given set and $\cH$ be a linear space of real functions defined on $\Omega$ such that if
$X_1,\dots,X_n\in\cH$, then $\phi(X_1,\dots,X_n)\in\cH$ for each $\phi\in lip(\mR^n)$. If $X\in\cH$, we call $X$ as a random variable.

\bd\label{subexp}
If a functional $\bar{\mE}: \cH\mapsto \mR$ satisfies: for $X, Y\in\cH$,

(i) $X\geq Y, \bar{\mE}[X]\geq\bar{\mE}[Y]$,

(ii) $\bar{\mE}[X+Y]\leq\bar{\mE}[X]+\bar{\mE}[Y]$,

(iii)$ \mbox{ for all}~\lambda\geq 0, \bar{\mE}[\lambda X]=\lambda\bar{\mE}[X]$,

(iv)$ \mbox{ for all}~c\in\mR, \bar{\mE}[X+c]=\bar{\mE}[X]+c$,

we call $\bar{\mE}$ as a sublinear expectation on $\cH$ and $(\Omega, \cH, \bar{\mE})$ as a sublinear expectation space.
\ed

Next, we define the distribution of a random vector on $(\Omega, \cH, \bar{\mE})$. For a $n$-dimensional random vector $X=(X_1, X_2, \cdots, X_n)$ for $X_i\in\cH$, $i=1,2, \cdots, n$, set
$$
F_X(\phi):=\bar{\mE}(\phi(X)), \qquad \phi\in lip(\mR^n),
$$
and then we call $F_X$ as the distribution of $X$.

\bd\label{samedist}
Assume that $X_1, X_2$ are two $n$-dimensional random vectors defined on different sublinear expectation spaces. If for all $\phi\in lip(\mR^n)$,
$$
 F_{X_1}(\phi)=F_{X_2}(\phi), 
 $$
 we say that the distributions of $X_1, X_2$ are the same.
\ed

\bd\label{indepen}
For two random vectors $Y=(Y_1, Y_2, \cdots, Y_m)$ for $Y_j\in\cH$ and $X=(X_1, X_2, \cdots, X_n)$ for $X_i\in\cH$, if for all $\phi\in lip(\mR^n\times\mR^m)$,
$$
\bar{\mE}[\phi(X,Y)]=\bar{\mE}[\bar{\mE}[\phi(x,Y)]_{x=X}],
$$
we say that $Y$ is independent from $X$.
\ed

Here, we use two above concepts to define L\'evy processes on $(\Omega, \cH, \bar{\mE})$.

\bd\label{defle}
Let $X=(X_t)_{t\geq0}$ be a $d$-dimensional c\`adl\`ag process on $(\Omega, \cH, \bar{\mE})$. If $X$ satisfies

(i) $X_0=0$;

(ii) for $t, s\geq 0$, the increment $X_{s+t}-X_t$ is independent from $(X_{t_1}, X_{t_2}, \cdots, X_{t_n})$, for any $n$ and $0\leq t_1<t_2\cdots<t_n\leq t$;

(iii) the distribution of $X_{s+t}-X_t$ does not depend on $t$;

we call $X$ as a L\'evy process.
\ed

\bd\label{defgle}
Assume that $X$ is a $d$-dimensional L\'evy process. If there exists a decomposition $X_t=X_t^c+X_t^d$ for $t\geq 0$, where $(X_t^c, X_t^d)$ is a $2d$-dimensional L\'evy process satisfying 
$$
\lim\limits_{t\downarrow 0}\frac{\bar{\mE}|X_t^c|^3}{t}=0, \quad \bar{\mE}|X_t^d|\leq Ct, \quad t\geq 0, \quad C\geq 0,
$$
we call $X$ as a G-L\'evy process.
\ed

In the following, we characterize G-L\'evy processes by partial differential equations.

\bt\label{gleop}
Assume that $X$ is a $d$-dimensional G-L\'evy process. Then for $g\in C_b^3(\mR^d)$ with $g(0)=0$, set
$$
G_X[g(\cdot)]:=\lim\limits_{t\downarrow 0}\frac{\bar{\mE}[g(X_t)]}{t},
$$
and then $G_X$ has the following L\'evy-Khintchine representation
\be
G_X[g(\cdot)]=\sup\limits_{(\nu,\zeta,Q)\in\cU}\left\{\int_{\mR^d\setminus\{0\}}g(u)\nu(\dif u)+\<\partial_x g(0), \zeta\>+\frac{1}{2}tr[\partial_x^2g(0)QQ^*]\right\},
\label{levykhin}
\ee
where $\cU$ is a subset of $\cM(\mR^d\setminus\{0\})\times\mR^d\times\mR^{d\times d}$, $\cM(\mR^d\setminus\{0\})$ is the collection of all measures on $(\mR^d\setminus\{0\}, \sB(\mR^d\setminus\{0\}))$, $\mR^{d\times d}$ is the set of all $d\times d$ matrices and $\cU$ satisfies 
\be
 \sup\limits_{(\nu,\zeta,Q)\in\cU}\left\{\int_{\mR^d\setminus\{0\}}|u|\nu(\dif u)+|\zeta|+\frac{1}{2}tr[QQ^*]\right\}<\infty.
 \label{usatcon}
 \ee
\et

\bt\label{glevis}
Suppose that $X$ is a $d$-dimensional G-L\'evy process. Then for $\phi\in C_{b, lip}(\mR^d)$, $v(t,x):=\bar{\mE}[\phi(x+X_t)]$ is the unique viscosity solution of the following partial integro-differential equation:
\ce
0&=&\partial_t v(t,x)-G_X[v(t,x+\cdot)-v(t,x)]\\
&=&\partial_t v(t,x)-\sup\limits_{(\nu,\zeta,Q)\in\cU}\bigg\{\int_{\mR^d\setminus\{0\}}[v(t,x+u)-v(t,x)]\nu(\dif u)+\<\partial_x v(t,x), \zeta\>\\
&&\qquad\qquad\qquad\qquad+\frac{1}{2}tr[\partial_x^2v(t,x)QQ^*]\bigg\}
\de
with the initial condition $v(0,x)=\phi(x)$.
\et

Conversely, if we have a set $\cU$ satisfying (\ref{usatcon}), is there a $d$-dimensional G-L\'evy process having the L\'evy-Khintchine representation (\ref{levykhin}) with the same set $\cU$? The answer is affirmed. We take $\Omega:=D_0(\mR^+, \mR^d)$, where $D_0(\mR^+, \mR^d)$ is the space of all c\`adl\`ag functions $\mR_+\ni t\mapsto \omega_t\in\mR^d$ with $\omega_0=0$, equipped with the Skorokhod topology. 

\bt\label{exisgle}
Suppose that $\cU$ satisfies (\ref{usatcon}). Then there exists a sublinear expectation $\bar{\mE}$ on $\Omega$ such that the canonical process $X$ is a $d$-dimensional G-L\'evy process having the L\'evy-Khintchine representation (\ref{levykhin}) with the same set $\cU$.
\et

\subsection{A capacity}\label{capa}

In the subsection, we introduce a capacity and related definitions.

First of all, fix a set $\cU$ satisfying (\ref{usatcon}) and $T>0$ and take $\Omega_T:=D_0([0,T], \mR^d)$ and the sublinear expectation $\bar{\mE}$ in Theorem \ref{exisgle}. Thus, we know that $(\Omega_T, \cH, \bar{\mE})$ is a sublinear expectation space. Here, we work on the space. Let $L_G^p(\Omega_T)$ be the completion of $lip(\Omega_T)$ under the norm $\|\cdot\|_p:=(\bar{\mE}|\cdot|^p)^{1/p}, p\geq 1$.

Let 
$$
\cV:=\{\nu\in\cM(\mR^d\setminus\{0\}): \exists (\zeta, Q)\in\mR^d\times\mR^{d\times d} ~\mbox{such ~ that}~ (\nu, \zeta, Q)\in\cU\}
$$
and let $\cG$ be the set of all the Borel measurable functions $g: \mR^d\mapsto\mR^d$ with $g(0)=0$.

{\bf Assumption:}
\begin{enumerate}[(i)]
\item
There exists a measure $\mu\in\cM(\mR^d)$ such that 
\ce
\int_{\mR^d\setminus\{0\}}|z|\mu(\dif z)<\infty, \quad \mu(\{0\})=0，
\de
and for all $\nu\in\cV$ there exists a function $g_{\nu}\in\cG$ satisfying 
$$
\nu(A)=\mu(g_{\nu}^{-1}(A)), \quad \forall A\in\sB(\mR^d\setminus\{0\}).
$$
\item
There exists $0<q<1$ such that 
$$
\sup\limits_{\nu\in\cV}\int_{0<|z|<1}|z|^q\nu(\dif z)<\infty.
$$
\item
$$
\sup\limits_{\nu\in\cV}\nu(\mR^d\setminus\{0\})<\infty.
$$
\end{enumerate}

Let $(\tilde{\Omega}, \sF, \mP)$ be a probability space supporting a Brownian motion $W$ and a Poisson random measure $N(\dif t, \dif z)$ with the intensity measure $\mu(\dif z)\dif t$. Let 
$$
\sF_t:=\sigma\left\{W_s, N((0,s],A): 0\leq s\leq t, A\in\sB(\mR^d\setminus\{0\})\right\}\vee\cN, \quad \cN:=\{U\in\sF, \mP(U)=0\}.
$$
We introduce the following set.

\bd\label{autt}
$\cA_{0,T}^{\cU}$ is a set of all the processes $\t_t=(\t^d_t, \t_t^{1,c},  \t_t^{2,c})$ for $t\in[0,T]$ satisfying 

(i) $(\t_t^{1,c},  \t_t^{2,c})$ is a $\sF_t$-adapted process and $\t^d$ is a $\sF_t$-predictable random field on $[0,T]\times\mR^d$,

(ii) For $\mP$-a.s. $\omega$ and a.e. $t\in[0,T]$,
$$
(\t^d(t,\cdot)(\omega), \t^{1,c}_t(\omega), \t^{2,c}_t(\omega))\in\left\{(g_{\nu},\zeta,Q)\in\cG\times\mR^d\times\mR^{d\times d}: (\nu,\zeta,Q)\in\cU\right\},
$$

(iii) 
$$
\mE^{\mP}\left[\int_0^T\(|\t^{1,c}_t|+\|\t^{2,c}_t\|^2+\int_{\mR^d\setminus\{0\}}|\t^d(t,z)|\mu(\dif z)\)\dif t\right]<\infty.
$$
\ed

For $\t\in\cA_{0,T}^{\cU}$, set
$$
B_t^{0,\t}:=\int_0^t\t^{1,c}_s\dif s+\int_0^t\t^{2,c}_s\dif W_s+\int_0^t\int_{\mR^d\setminus\{0\}}\t^d(s,z)N(\dif s,\dif z), \quad t\in[0,T],
$$
and by Corollary 14 in \cite{pk1}, it holds that for $\xi\in L_G^1(\Omega_T)$
$$
\bar{\mE}[\xi]=\sup\limits_{\t\in\cA_{0,T}^{\cU}}\mE^{\mP^{\t}}[\xi], \quad \mP^{\t}=\mP\circ (B_{\cdot}^{0,\t})^{-1}.
$$
And then define 
$$
\bar{C}(D):=\sup\limits_{\t\in\cA_{0,T}^{\cU}}\mP^{\t}(D), \quad D\in\sB(\Omega_T),
$$
and $\bar{C}$ is a capacity. For $D\in\sB(\Omega_T)$, if $\bar{C}(D)=0$, we call $D$ as a polar set. So, if a property holds outside a polar set, we say that the property holds quasi-surely (q.s. in short).

\subsection{The It\^o integrals with respect to G-L\'evy processes}\label{itointe}

In the subsection, we introduce the It\^o integrals with respect to G-L\'evy processes under the framework of the above subsection. 

Let $X$ denote the canonical process on the space, i.e. $X_t(\omega)=\omega_t, t\in[0,T]$. So, $X$ is a $d$-dimensional G-L\'evy process. Although the It\^o integrals with respect to G-Brownian motions have been introduced in \cite{q1}, we need to introduce two related spaces used in the sequel. Take $0=t_0<t_1<\cdots<t_N=T$. Let $p\geq 1$ be fixed.
Define 
\ce
\cM^{p,0}_G(0,T):=\Big\{\eta_t(\omega)=\sum\limits_{j=1}^N\xi_{j-1}(\omega)1_{[t_{j-1},t_j)}(t);\xi_{j-1}(\omega)\in L^p_G(\Omega_{t_{j-1}})\Big\}.
\de 
Let $\cM^p_G(0,T)$ and $\cH_G^p(0,T)$ denote the completion of $\cM^{p,0}_G(0,T)$ under the norm
$$
\|\eta\|_{\cM^p_G(0,T)}=\left(\int_0^T\bar{\mE}|\eta_t|^p\dif t\right)^{\frac{1}{p}} ~\mbox{and}~\|\eta\|_{\cH^p_G(0,T)}=\left(\bar{\mE}\left(\int_0^T|\eta_t|^2\dif t\right)^{\frac{p}{2}}\right)^{\frac{1}{p}},
$$
respectively. Let $\cM^p_G([0,T], \mR^d)$ and $\cH^p_G([0,T], \mR^d)$ be the collection of all the processes 
$$
\eta_t=(\eta^1_t, \eta^2_t, \cdots, \eta^d_t), \quad t\in[0,T], \quad \eta^i\in \cM^{p}_G(0,T) ~\mbox{and}~ \cH_G^p(0,T),
$$
respectively. 

Next, we introduce the It\^o integrals with respect to random measures. First, define a random measure: for any $0\leq t\leq T$ and $A\in\sB(\mR^d\setminus\{0\})$,
\ce
\kappa_t:=X_t-X_{t-}, \quad L((0,t], A):=\sum\limits_{0<s\leq t}I_A(\kappa_s), \quad q.s.. 
\de
And then we define the It\^o integral with respect to the random measure $L(\dif t, \dif u)$. Let $\cH^S_G([0,T]\times(\mR^d\setminus\{0\}))$ be the collection of all the processes defined on $[0,T]\times(\mR^d\setminus\{0\})\times\Omega$ with the form
$$
f(s,u)(\omega)=\sum\limits_{k=1}^{n-1}\sum\limits_{l=1}^m\phi_{k,l}(X_{t_1}, X_{t_2}-X_{t_1}, \cdots, X_{t_k}-X_{t_{k-1}})I_{[t_k, t_{k+1})}(s)\psi_l(u), n,m\in\mN,
$$
where $0\leq t_1<\cdots<t_n\leq T$ is a partition of $[0,T]$, $\phi_{k,l}\in C_{b, lip}(\mR^{d\times k})$ and $\{\psi_l\}_{l=1}^m\subset C_{b, lip}(\mR^d)$ are functions with disjoint supports and $\psi_l(0)=0$.

\bd\label{simito}
For any $f\in\cH^S_G([0,T]\times(\mR^d\setminus\{0\}))$, set 
$$
\int_0^t\int_{\mR^d\setminus\{0\}}f(s,u)L(\dif s, \dif u):=\sum\limits_{0<s\leq t}f(s,\kappa_s)I_{\mR^d\setminus\{0\}}(\kappa_s), \quad q.s..
$$ 
\ed

By Theorem 28 in \cite{pk1}, we have that $\int_0^t\int_{\mR^d\setminus\{0\}}f(s,u)L(\dif s, \dif u)\in L_G^2(\Omega_T)$. Let $\cH^2_G([0,T]\times(\mR^d\setminus\{0\}))$ be the completion of $\cH^S_G([0,T]\times(\mR^d\setminus\{0\}))$ with respect to the norm $\|\cdot\|_{\cH^2_G([0,T]\times(\mR^d\setminus\{0\}))}$, where 
$$
\|f\|_{\cH^2_G([0,T]\times(\mR^d\setminus\{0\}))}:=\bar{\mE}\left[\int_0^T\sup\limits_{\nu\in\cV}\int_{\mR^d\setminus\{0\}}|f(s,u)|^2\nu(\dif u)\dif s\right]^{1/2}, \quad f\in\cH^S_G([0,T]\times(\mR^d\setminus\{0\})).
$$
Thus, Corollary 29 in \cite{pk1} admits us to get that for $f\in\cH^2_G([0,T]\times(\mR^d\setminus\{0\}))$, 
\be
\int_0^t\int_{\mR^d\setminus\{0\}}f(s,u)L(\dif s, \dif u)=\sum\limits_{0<s\leq t}f(s,\kappa_s)I_{\mR^d\setminus\{0\}}(\kappa_s), \quad q.s..
\label{jumpin}
\ee
Let $\cH^2_G([0,T]\times(\mR^d\setminus\{0\}), \mR^d)$ be the space of all the processes 
$$
f(t,u)=\(f^1(t,u), f^2(t,u), \cdots, f^d(t,u)\), \quad f^i\in\cH^2_G([0,T]\times(\mR^d\setminus\{0\})).
$$

\subsection{Stochastic differential equations driven by G-L\'evy processes}\label{sdegle}

In the subsection, we introduce stochastic differential equations driven by G-L\'evy processes and related additive functionals. 

First, we introduce some notations. Let $\mS^d$ be the space of all $d\times d$ symmetric matrices. For $A\in\mS^d$, set
$$
G(A):=\frac{1}{2}\sup\limits_{Q\in\cQ}\tr[QQ^*A],
$$
where $\cQ$ is a nonempty, bounded, closed and convex subset of $\mR^{d\times d}$. And then $G: \mS^d\mapsto\mR$ is a monotonic, sublinear and positive homogeneous functional(\cite{peng}). We choose $\cU\subset\cM(\mR^d\setminus\{0\})\times\{0\}\times\cQ$ satisfying (\ref{usatcon}) and still work under the framework of Subsection \ref{capa}. So, the canonical process $X_t$ can be represented as $X_t=B_t+X_t^d$, where $B_t$ is a G-Brownian motion associated with $\cQ$ and $X_t^d$ is a pure jump G-L\'evy process associated with $\cM(\mR^d\setminus\{0\})$.

Next, we consider Eq.(\ref{glesde}) and assume:
\begin{enumerate}[(${\bf H}^1_{b,h,\sigma,f}$)]
\item 
There exists a constant $C_1>0$ such that for any $t\in[0,T]$ and $x,y\in\mR^d$, 
\ce
&&|b(t,x)-b(t,y)|^2+|h_{ij}(t,x)-h_{ij}(t,y)|^2+\|\sigma(t,x)-\sigma(t,y)\|^2\\
&+&\sup\limits_{\nu\in\cV}\int_{\mR^d\setminus\{0\}}|f(t,x,u)-f(t,y,u)|^2\nu(\dif u)\leq C_1\rho(|x-y|^2),
\de
where $\rho:(0,+\infty)\mapsto(0,+\infty)$ is a continuous, increasing
and concave function so that
\ce
\rho(0+)=0, \quad
\int_0^1\frac{dr}{\rho(r)}=+\infty.
\de
\end{enumerate}
\begin{enumerate}[(${\bf H}^2_{b,h,\sigma,f}$)]
\item 
There exists a constant $C_2>0$ such that for any $t\in[0,T]$
\ce
|b(t,0)|^2+|h_{ij}(t,0)|^2+\|\sigma(t,0)\|^2+\sup\limits_{\nu\in\cV}\int_{\mR^d\setminus\{0\}}|f(t,0,u)|^2\nu(\dif u)\leq C_2.
\de
\end{enumerate}

By \cite[Theorem 3.1]{wg}, we know that under (${\bf H}^1_{b,h,\sigma,f}$)-(${\bf H}^2_{b,h,\sigma,f}$), Eq.(\ref{glesde}) has a unique solution $Y_t$ with 
\be
\bar{\mE}\left[\sup\limits_{t\in[0,T]}|Y_t|^2\right]<\infty.
\label{estiy}
\ee

And then we introduce the following additive functional
\be
F_{s,t}&:=&\a\int_s^t G(g_1)(r,Y_r)\dif r+\beta\int_s^t g^{ij}_1(r,Y_r)\dif \<B^i,B^j\>_r+\int_s^t \<g_2(r,Y_r),\dif B_r\>\no\\
&&+\int_s^t\int_{\mR^d\setminus\{0\}}g_3(r,Y_r, u)L(\dif r, \dif u)+\int_s^t\sup\limits_{\nu\in\cV}\int_{\mR^d\setminus\{0\}}\gamma g_3(r,Y_r, u)\nu(\dif u)\dif r, \no\\
&&\quad 0\leq s<t\leq T,
\label{addfun}
\ee
where $\a, \beta, \gamma\in\mR$ are three constants and
\ce
&&g_1: [0,T]\times\mR^d\mapsto\mR^{d\times d}, \quad g_1^{ij}=g_1^{ji}, \\
&&g_2:[0,T]\times\mR^d\mapsto\mR^d, \\
&&g_3:[0,T]\times\mR^d\times(\mR^d\setminus\{0\})\mapsto\mR,
\de
are Borel measurable so that $F_{s,t}$ is well-defined.

\bd\label{pathinde}
The additive functional $F_{s,t}$ is called path independent, if there exists a function 
$$
V: [0,T]\times\mR^d\mapsto\mR,
$$
such that for any $s\in[0, T]$ and $Y_s\in L^2_G(\Omega_T)$, the solution $(Y_t)_{t\in[s,T]}$ of Eq.(\ref{glesde}) satisfies
\be
F_{s,t}=V(t,Y_t)-V(s,Y_s).
\label{defi}
\ee
\ed

\section{Main results and their proofs}\label{main}

In the section, we state and prove the main results under the framework of Subsection \ref{sdegle}. And then we analysis some special cases and compare our result with some known results.

\subsection{Main results}\label{maires} 

In the subsection, we state and prove the main results. Let us begin with a key lemma.

\bl\label{zero}
Assume that $\cQ$ is bounded away from $0$ and $Z_t$ is a $1$-dimensional G-It\^o-L\'evy process, i.e.
\be
Z_t=\int_0^t\Gamma_s\dif s+\int_0^t\Phi_{ij}(s)\dif \<B^i, B^j\>_s+\int_0^t\<\Psi_s,\dif B_s\>+\int_0^t\int_{\mR^d\setminus\{0\}}K(s,u)L(\dif s, \dif u),
\label{gitole}
\ee
where $\Gamma\in\cM^1_G(0,T), \Phi_{ij}\in\cM^1_G(0,T), \Phi_{ij}=\Phi_{ji}, i,j=1,2,\cdots,d, \Psi\in\cH^1_G([0,T], \mR^d), K\in\cH^2_G([0,T]\times(\mR^d\setminus\{0\}))$.
Then $Z_t=0$ for all $t\in[0,T]$ q.s. if and only if $\Gamma_t=0, \Phi_{ij}(t)=0, \Psi_t=0$ a.e.$\times$q.s. on $[0,T]\times\Omega_T$ and $K(t,u)=0$ a.e.$\times$a.e.$\times$q.s. on $[0,T]\times(\mR^d\setminus\{0\})\times\Omega_T$.
\el
\begin{proof}
Sufficiency is direct if one inserts $\Gamma_t=0, \Phi_{ij}(t)=0, \Psi_t=0, K(t,u)=0$ into (\ref{gitole}). Let us prove necessity. If $Z_t=0$ for any $t\in[0,T]$, we get that
\be
0=\int_0^t\Gamma_s\dif s+\int_0^t\Phi_{ij}(s)\dif \<B^i, B^j\>_s+\int_0^t\Psi_s^*\dif B_s+\int_0^t\int_{\mR^d\setminus\{0\}}K(s,u)L(\dif s, \dif u).
\label{gitole1}
\ee
By taking the quadratic process with $\int_0^t\<\Psi_s, \dif B_s\>$ on two sides of (\ref{gitole1}), it holds that
\ce
0&=&\<\int_0^{\cdot}\Psi^*_s\dif B_s, \int_0^{\cdot}\Psi^*_s\dif B_s\>_t=\int_0^t\Psi^{i}_s\Psi^{j}_s\dif \<B^i, B^j\>_s=\int_0^t\tr(\Psi_s\Psi^*_s\dif \<B\>_s)\\
&=&\int_0^t\tr(\dif \<B\>_s\Psi_s\Psi^*_s)=\int_0^t\<\dif \<B\>_s\Psi_s, \Psi_s\>,
\de
where 
\ce
\<B\>:=\left(\begin{array}{c}
\<B^1,B^1\>\quad \<B^1,B^2\>\cdots\<B^1,B^d\>\\
\vdots\qquad\qquad\vdots\qquad\qquad\vdots\\
\<B^d,B^1\>\quad \<B^d,B^2\>\cdots\<B^d,B^d\>
\end{array}
\right).
\de
Note that $\cQ$ is bounded away from $0$. Thus, there exists a constant $\iota>0$ such that $\<B\>_s\geq\iota sI_d$ and then 
$$
0=\int_0^t\<\dif \<B\>_s\Psi_s, \Psi_s\>\geq \iota\int_0^t\<\Psi_s, \Psi_s\>\dif s.
$$
From this, we know that $\Psi_t=0$ a.e.$\times$q.s.. So, (\ref{gitole1}) becomes
$$
0=\int_0^t\Gamma_s\dif s+\int_0^t\Phi_{ij}(s)\dif \<B^i, B^j\>_s+\int_0^t\int_{\mR^d\setminus\{0\}}K(s,u)L(\dif s, \dif u).
$$

Next, set
$$
\tau_0=0, \quad \tau_n:=\inf\{t>\tau_{n-1}: \kappa_t\neq 0 \}, \quad n=1,2,\cdots,
$$
and $\{\tau_n\}$ is a stopping time sequence with respect to $(\sB_t)_{t\geq0}$ and $\tau_n\uparrow\infty$ as $n\rightarrow\infty$ q.s.(c.f. \cite[Proposition 16]{pk1}). So, by (\ref{jumpin}) it holds that for $t\in[0, \tau_1\land T)$, 
\ce
0=\int_0^{\tau_1\land T}\Gamma_s\dif s+\int_0^{\tau_1\land T}\Phi_{ij}(s)\dif \<B^i, B^j\>_s,
\de
i.e.
$$
-\int_0^{\tau_1\land T}\Gamma_s\dif s=\int_0^{\tau_1\land T}\Phi_{ij}(s)\dif \<B^i, B^j\>_s.
$$
Thus, by the similar deduction to that in \cite[Corollary 1]{ry} one can have that 
$$
\bar{\mE}\int_0^{\tau_1\land T}\left(\tr[\Phi_s\Phi_s]\right)^{1/2}\dif s=\bar{\mE}\int_0^{\tau_1\land T}|\Gamma_s|\dif s=0.
$$
Based on this, we know that $\Phi_t=0, \Gamma_t=0$ for $t\in[0, \tau_1\land T)$. If $\tau_1\geq T$, the proof is over; if $\tau_1<T$, we continue. For $t=\tau_1$, (\ref{gitole1}) goes to 
$$
\Phi_t=0, \quad \Gamma_t=0, \quad K(t,\kappa_t)=0.
$$
For $t\in[\tau_1, \tau_2\land T)$, by the same means to the above for $t\in[0, \tau_1\land T)$, we get that $\Phi_t=0, \Gamma_t=0$ for $t\in[\tau_1, \tau_2\land T)$. If $\tau_2\geq T$, the proof is over; if $\tau_2<T$, we continue till $T\leq\tau_n$. Thus, we obtain that $\Gamma_t=0, \Phi_{ij}(t)=0, \Psi_t=0$ a.e.$\times$q.s. on $[0,T]\times\Omega_T$ and $K(t,u)=0$ a.e.$\times$a.e.$\times$q.s. on $[0,T]\times(\mR^d\setminus\{0\})\times\Omega_T$. The proof is complete.
\end{proof}

The main result in the section is the following theorem.

\bt\label{suffnece}
Assume that $\cQ$ is bounded away from $0$ and $b, h, \sigma, f$ satisfy (${\bf H}^1_{b,h,\sigma,f}$)-(${\bf H}^2_{b,h,\sigma,f}$). Then for $V\in C^{1,2}_b([0,T]\times\mR^d)$, $F_{s,t}$ is path independent in the sense of (\ref{defi}) if and only if $(V, g_1, g_2, g_3)$ satisfies the partial integral-differential equation
\be\left\{\begin{array}{ll}
\partial_tV(t,x)+\<\partial_xV(t,x), b(t,x)\>=\a G(g_1)(t,x)+\sup\limits_{\nu\in\cV}\int_{\mR^d\setminus\{0\}}\gamma g_3(t,x, u)\nu(\dif u),\\
\<\partial_xV(t,x), h_{ij}(t,x)\>+\frac{1}{2}\<\partial^2_xV(t,x)\sigma^i(t,x), \sigma^j(t,x)\>=\beta g^{ij}_1(t, x),\\
(\sigma^T\partial_x V)(t, x)=g_2(t, x),\\
V\(t, x+f(t,x,u)\)-V(t, x)=g_3(t, x,u), \\
t\in[0,T], x\in\mR^d, u\in\mR^d\setminus\{0\}.
\end{array}
\label{three}
\right.
\ee
\et
\begin{proof}
First, we prove necessity. On one hand, since $F_{s,t}$ is path independent in the sense of (\ref{defi}), by Definition \ref{pathinde} it holds that
\be
V(t,Y_t)-V(s,Y_s)&=&\a\int_s^t G(g_1)(r,Y_r)\dif r+\beta\int_s^t g^{ij}_1(r,Y_r)\dif \<B^i,B^j\>_r+\int_s^t \<g_2(r,Y_r),\dif B_r\>\no\\
&&+\int_s^t\int_{\mR^d\setminus\{0\}}g_3(r,Y_r, u)L(\dif r, \dif u)+\int_s^t\sup\limits_{\nu\in\cV}\int_{\mR^d\setminus\{0\}}\gamma g_3(r,Y_r, u)\nu(\dif u)\dif r.\no\\
\label{pathindepen}
\ee
On the other hand, applying the It\^o formula for G-It\^o-L\'evy processes (\cite[Theorem 32]{pk1}) to $V(t,Y_t)$, one can obtain that
\be
V(t,Y_t)-V(s,Y_s)&=&\int_s^t\partial_rV(r,Y_r)\dif r+\int_s^t\partial_kV(r,Y_r)b^k(r,Y_r)\dif r\no\\
&&+\int_s^t\partial_kV(r,Y_r)h^k_{ij}(r,Y_r)\dif \<B^i, B^j\>_r+\int_s^t\<(\sigma^T\partial_xV)(r,Y_r),\dif B_r\>\no\\
&&+\int_s^t\int_{\mR^d\setminus\{0\}}\(V(r,Y_r+f(r,Y_r,u))-V(r,Y_r)\)L(\dif r,\dif u)\no\\
&&+\frac{1}{2}\int_s^t\partial_{kl}V(r,Y_r)\sigma^{ki}(r,Y_r)\sigma^{lj}(r,Y_r)\dif \<B^i, B^j\>_r.
\label{itofor}
\ee
By (\ref{estiy}) (${\bf H}^1_{b,h,\sigma,f}$)-(${\bf H}^2_{b,h,\sigma,f}$), one can verify that 
\ce
&&\partial_rV(r,Y_r)+\partial_kV(r,Y_r)b^k(r,Y_r)\in\cM^1_G(0,T),\\
&&\partial_kV(r,Y_r)h^k_{ij}(r,Y_r)+\frac{1}{2}\partial_{kl}V(r,Y_r)\sigma^{ki}(r,Y_r)\sigma^{lj}(r,Y_r)\in\cM^1_G(0,T),\\
&&(\sigma^T\partial_xV)(r,Y_r)\in\cH^1_G([0,T], \mR^d),\\
&&V(r,Y_r+f(r,Y_r,u))-V(r,Y_r)\in\cH^2_G([0,T]\times(\mR^d\setminus\{0\})).
\de
Thus, by (\ref{pathindepen}) (\ref{itofor}) and Lemma \ref{zero} we know that
\ce\left\{\begin{array}{ll}
\partial_rV(r,Y_r)+\<\partial_xV(r,Y_r), b(r,Y_r)\>=\a G(g_1)(r,Y_r)+\sup\limits_{\nu\in\cV}\int_{\mR^d\setminus\{0\}}\gamma g_3(r,Y_r, u)\nu(\dif u),\\
\<\partial_xV(r,Y_r), h_{ij}(r,Y_r)\>+\frac{1}{2}\<\partial^2_xV(r,Y_r)\sigma^i(r,Y_r), \sigma^j(r,Y_r)\>=\beta g^{ij}_1(r,Y_r),\\
(\sigma^T\partial_x V)(r,Y_r)=g_2(r,Y_r),\\
V\(r,Y_r+f(r,Y_r,u)\)-V(r,Y_r)=g_3(r,Y_r,u), \quad a.e.\times q.s..
\end{array}
\right.
\de
Now, we insert $r=s, Y_s=x\in\mR^d$ into the above equalities and get that
\ce\left\{\begin{array}{ll}
\partial_sV(s,x)+\<\partial_xV(s,x), b(s,x)\>=\a G(g_1)(s,x)+\sup\limits_{\nu\in\cV}\int_{\mR^d\setminus\{0\}}\gamma g_3(s,x, u)\nu(\dif u),\\
\<\partial_xV(s,x), h_{ij}(s,x)\>+\frac{1}{2}\<\partial^2_xV(s,x)\sigma^i(s,x), \sigma^j(s,x)\>=\beta g^{ij}_1(s, x),\\
(\sigma^T\partial_x V)(s, x)=g_2(s, x),\\
V\(s, x+f(s,x,u)\)-V(t, x)=g_3(s, x,u), \\
s\in[0,T], x\in\mR^d, u\in\mR^d\setminus\{0\}.
\end{array}
\right.
\de
Since $s, x$ are arbitrary, we have (\ref{three}).

Next, we treat sufficiency. By the It\^o formula for G-It\^o-L\'evy processes to $V(t,Y_t)$, we have (\ref{itofor}). And then one can apply (\ref{three}) to (\ref{itofor}) to get (\ref{pathindepen}). That is, $F_{s,t}$ is path independent in the sense of (\ref{defi}). The proof is complete.
\end{proof}

\subsection{Some special cases}

In the subsection, we analysis some special cases. 

If $b(t,x)=0, h_{ij}(t,x)=0, \sigma(t,x)=I_d, f(t,x,u)=u$, Eq.(\ref{glesde}) becomes 
\ce
\dif Y_t=\dif B_t+\int_{\mR^d\setminus\{0\}}uL(\dif t, \dif u)=\dif X_t.
\de
We take $\a=1, \beta=\frac{1}{2}, \gamma=1$. Thus, by Theorem \ref{suffnece}, it holds that $F_{s,t}$ is path independent in the sense of (\ref{defi}) if and only if $(V, g_1, g_2, g_3)$ satisfies the partial integral-differential equation
\ce\left\{\begin{array}{ll}
\partial_tV(t,x)=G(\partial^2_xV)(t,x)+\sup\limits_{\nu\in\cV} \int_{\mR^d\setminus\{0\}}(V(t, x+u)-V(t, x))\nu(\dif u),\\
\partial^2_xV(t,x)= g_1(t, x),\\
\partial_x V(t, x)=g_2(t, x),\\
V(t, x+u)-V(t, x)=g_3(t, x,u), \\
t\in[0,T], x\in\mR^d, u\in\mR^d\setminus\{0\}.
\end{array}
\right.
\de
Besides, by Theorem \ref{glevis}, it holds that for $\phi\in C_{b, lip}(\mR^d)$, $V(t,x)=\bar{\mE}[\phi(x+X_t)]$ is the unique viscosity solution of the following partial integro-differential equation:
\ce
\partial_t V(t,x)-G(\partial^2_xV)(t,x)-\sup\limits_{\nu\in\cV}\int_{\mR^d\setminus\{0\}}\(V(t, x+u)-V(t, x)\)\nu(\dif u)=0
\de
with the initial condition $V(0,x)=\phi(x)$. So, 
\ce
&&g_1(t, x)=\partial^2_x\bar{\mE}[\phi(x+X_t)], \\
&&g_2(t, x)=\partial_x\bar{\mE}[\phi(x+X_t)], \\
&&g_3(t, x,u)=\bar{\mE}[\phi(x+u+X_t)]-\bar{\mE}[\phi(x+X_t)].
\de
That is, we can find $g_1, g_2, g_3$. 

If $d=1, \a=1, \beta=0, \gamma=0$, it follows from Theorem \ref{suffnece} that 
$F_{s,t}$ is path independent in the sense of (\ref{defi}) if and only if $(V, g_1, g_2, g_3)$ satisfies the partial integral-differential equation
\ce\left\{\begin{array}{ll}
\partial_tV(t,x)+\partial_xV(t,x)b(t,x)=G(g_1)(t,x),\\
\partial_xV(t,x)h(t,x)+\frac{1}{2}\partial^2_xV(t,x)\sigma^2(t,x)=0,\\
(\sigma\partial_x V)(t, x)=g_2(t, x),\\
V\(t, x+f(t,x,u)\)-V(t, x)=g_3(t, x,u), \\
t\in[0,T], x\in\mR, u\in\mR\setminus\{0\}.
\end{array}
\right.
\de
And then if $\sigma(t,x)\neq 0$, the unique solution of the above second equation is 
$$
V(t,x)=V(t,0)+\partial_xV(t,0)\int_0^xe^{-2\int_0^z\frac{h(t,v)}{\sigma^2(t,v)}\dif v}\dif z.
$$
Besides, since $\cQ$ is bounded away from $0$, $G$ is invertible. Thus, 
\ce
&&g_1(t,x)=G^{-1}\(\partial_tV(t,x)+b(t,x)\partial_xV(t,0)e^{-2\int_0^x\frac{h(t,v)}{\sigma^2(t,v)}\dif v}\),\\
&&g_2(t, x)=\sigma(t, x)\partial_xV(t,0)e^{-2\int_0^x\frac{h(t,v)}{\sigma^2(t,v)}\dif v},\\
&&g_3(t, x,u)=\partial_xV(t,0)\int_x^{x+f(t,x,u)}e^{-2\int_0^z\frac{h(t,v)}{\sigma^2(t,v)}\dif v}\dif z.
\de
That is, we also give out $g_1, g_2, g_3$ in the case. 

\br
In the above special cases, we can describe concretely $g_1, g_2, g_3$. This is interesting and also is one of our motivations. 
\er

\subsection{Comparison with some known results}\label{com}

In the subsection, we compare our result with some known results. 

First, if we take $f(t,x,u)=0$ in Eq.(\ref{glesde}) and $g_3(t,x,u)=0$ in (\ref{addfun}), Theorem \ref{suffnece} becomes \cite[Theorem 2]{ry}. Therefore, our result is more general.

Second, we take $\cM(\mR^d\setminus\{0\})=\{\nu\}, \cQ=\{I_d\}$ in Subsection \ref{sdegle}. Thus, $B$ is a classical Brownian motion with $\<B^i, B^j\>_t=I_{i=j}t$ and $L(\dif t, \dif u)$ is a classical Poisson random measure. And then Eq.(\ref{glesde}) goes into
\be
\dif Y_t=b(t,Y_t)\dif t +\sum\limits_{i=1}^d h_{ii}(t,Y_t)\dif t+\sigma(t,Y_t)\dif B_t+\int_{\mR^d\setminus\{0\}}f(t,Y_t,u)L(\dif t,\dif u).
\label{special}
\ee
Note that $\nu(\mR^d\setminus\{0\})<\infty$. So, by (\ref{estiy}) (${\bf H}^1_{b,h,\sigma,f}$)-(${\bf H}^2_{b,h,\sigma,f}$), \cite[Theorem 13]{pk2} admits us to obtain that 
$$
\int_0^t\int_{\mR^d\setminus\{0\}}f(s,Y_s,u)\tilde{L}(\dif s,\dif u):=\int_0^t\int_{\mR^d\setminus\{0\}}f(s,Y_s,u)L(\dif s,\dif u)-\int_0^t\int_{\mR^d\setminus\{0\}}f(s,Y_s,u)\nu(\dif u)\dif s
$$
is a $\sB_t$-martingale, where $\sB_t:=\sigma\{\omega_s, 0\leq s\leq t\}, 0\leq t\leq T$. Therefore, we can rewrite Eq.(\ref{special}) to get that
\ce
\dif Y_t&=&b(t,Y_t)\dif t +\sum\limits_{i=1}^d h_{ii}(t,Y_t)\dif t+\int_{\mR^d\setminus\{0\}}f(t,Y_t,u)\nu(\dif u)\dif t+\sigma(t,Y_t)\dif B_t\\
&&+\int_{\mR^d\setminus\{0\}}f(t,Y_t,u)\tilde{L}(\dif t,\dif u).
\de
This is a classical stochastic differential equation with jumps. Thus, by \cite[Theorem 1.2]{q0} it holds that under (${\bf H}^1_{b,h,\sigma,f}$)-(${\bf H}^2_{b,h,\sigma,f}$), the above equation has a unique solution.

In the following, note that $G(g_1)=\frac{1}{2}\tr(g_1)=\frac{1}{2}\sum\limits_{i=1}^d g_1^{ii}$. And then $F_{s,t}$ can be represented as
\ce
F_{s,t}&=&\int_s^t \left(\frac{\a}{2}+\beta\right)\sum\limits_{i=1}^d g_1^{ii}(r,Y_r)\dif r+\int_s^t \<g_2(r,Y_r),\dif B_r\>+\int_s^t\int_{\mR^d\setminus\{0\}}g_3(r,Y_r, u)\tilde{L}(\dif r, \dif u)\\
&&+\int_s^t\int_{\mR^d\setminus\{0\}}(1+\gamma )g_3(r,Y_r, u)\nu(\dif u)\dif r.
\de
This is just right \cite[(3)]{qw3} without the distribution of $Y_r$ for $r\in[s,t]$. So, in the case Definition \ref{pathinde} and Theorem \ref{suffnece} are Definition 2.1 and Theorem 3.2 in \cite{qw3} without the distribution of $Y_r$ for $r\in[s,t]$, respectively. 
Therefore, our result overlaps \cite[Theorem 3.2]{qw3} in some sense.

\bigskip

\textbf{Acknowledgements:}

The authors are very grateful to Professor Xicheng Zhang for valuable discussions. The first author also thanks Professor Renming Song for providing her an excellent environment to work in the University of Illinois at Urbana-Champaign.

\end{document}